\documentclass{tran-l}

\usepackage{graphicx}

\usepackage[cmtip,all]{xy}

\newtheorem{lemma}{Lemma}[section]
\newtheorem{theorem}[lemma]{Theorem}
\newtheorem{proposition}[lemma]{Proposition}

\theoremstyle{definition}
\newtheorem{definition}[lemma]{Definition}
\theoremstyle{remark}
\newtheorem{remark}[lemma]{Remark}

\numberwithin{equation}{section}

\usepackage[noabbrev]{cleveref}

\crefname{section}{Section}{Sections}
\crefname{lemma}{Lemma}{Lemmas}
\crefname{proposition}{Proposition}{Propositions}
\crefname{theorem}{Theorem}{Theorems}
\crefname{corollary}{Corollary}{Corollaries}
\crefname{definition}{Definition}{Definitions}
\crefname{equation}{Equation}{Equations}
\crefname{figure}{Figure}{Figures}

\newcommand{\iso}{\simeq}
\newcommand{\homot}{\approx}
\newcommand{\const}[1]{{\bullet_{#1}}}
\newcommand{\concat}{*}
\newcommand{\rev}[1]{\overline{#1}}

\newcommand{\loops}[2][]{\Omega^{#1}_{#2}}
\DeclareMathOperator{\length}{length}
\newcommand{\crit}{\lambda}

\DeclareMathOperator{\im}{im}

\DeclareMathOperator{\diam}{diam}
\DeclareMathOperator{\vol}{Vol}
\DeclareMathOperator{\supp}{supp}
\DeclareMathOperator{\inj}{inj}

\newcommand{\Q}{\mathbb{Q}}

\newcommand{\Sp}[1][2]{\mathbb{S}^{#1}}

\newcommand{\card}[1]{\left\lvert#1\right\rvert}
\newcommand{\dist}[2][]{\operatorname{dist}_{#1}(#2)}

\newcommand{\homto}[1]{\xrightarrow{#1}}

\begin{document}

\title[Linear length bounds on geodesics in closed Riemannian surfaces]{ Curvature-free linear length bounds on geodesics in closed Riemannian surfaces}


\author{Herng Yi Cheng}
\address{Department of Mathematics, University of Toronto, 40 St.~George Street, Toronto, Ontario, Canada M5S 2E4}
\curraddr{}
\email{herngyi@math.toronto.edu}
\thanks{This research has been partially supported by the Vivekananda Graduate Scholarship from the Department of Mathematics at the University of Toronto.}


\date{}

\dedicatory{}

\begin{abstract}
    This paper proves that in any closed Riemannian surface $M$ with diameter $d$, the length of the $k^\text{th}$-shortest geodesic between two given points $p$ and $q$ is at most $8kd$. This bound can be tightened further to $6kd$ if $p = q$. This improves prior estimates by A.~Nabutovsky and R.~Rotman.
\end{abstract}

\maketitle


\section{Introduction}
\label{sec:Introduction}

It was proven in the 1950s by J.-P.~Serre that any two points $p$ and $q$ on a closed Riemannian manifold $M$ are connected by infinitely many distinct geodesics \cite{Serre}. A natural question is to understand the growth rate of the lengths of the geodesics from $p$ to $q$. In particular, what bounds on the growth rate are dictated by topology, independent of geometric invariants such as curvature? (For the bounds to be scale-invariant, they could depend linearly on the diameter $d = \diam(M)$ or $\vol(M)^{1/n}$.)

In 1958, A.~Schwarz suggested a way to exploit the topology of $M$---or rather, the topology of $\loops{pq}M$, the space of piecewise smooth paths in $M$ from $p$ to $q$---to prove absolute bounds in his ``quantitative'' proof of Serre's theorem.\footnote{For a detailed definition and treatment of $\loops{pq}M$, see Part~III in \cite{MorseTheory_Milnor}.} He represented nontrivial homology classes of $\loops{pq}(M)$ using only paths of controlled length, and used them to construct at least $k$ non-constant geodesics from $p$ to $q$ with length at most $c(M)k$, where $c(M)$ is some constant depending on the metric on $M$ \cite{Schwartz}. In 2013, A.~Nabutovsky and R.~Rotman proved that there are at least $k$ geodesics from $p$ to $q$ of length at most $4nk^2d$, where $n = \dim(M)$ \cite{Nabutovsky_QuantitativeMorse}. Thus they eliminated the bound's dependence on geometric invariants other than the diameter, at the cost of a quadratic growth in $k$.

On the other hand, when $n = 2$, they managed to prove a linear bound on the growth rate. Their key results are linear bounds for Riemannian 2-spheres $M = (\Sp,g)$: there are at least $k$ geodesics from $p$ to $q$ of length at most $22(k - 1)d + \dist{p,q} \leq 22kd$, where $d = \diam(M)$ \cite{LinearBoundsSegments}. They also showed that $M$ has at least $k$ distinct non-constant geodesic loops based at $p$ of length at most $20kd$ \cite{LinearBounds}.

(A geodesic $\gamma : I \to M$ is called a \emph{geodesic loop} if $\gamma(0) = \gamma(1)$. In this case, if $\gamma'(0) = \gamma'(1)$ then $\gamma$ is also called a \emph{periodic geodesic}. If $\gamma(0) \neq \gamma(1)$ then $\gamma$ is called a \emph{geodesic segment}. Note that a geodesic segment may not be a shortest geodesic between $p$ and $q$.)

In a sense, the simply-connected $\Sp$ is the ``hardest case'', as an analysis of the fundamental group of a closed Riemannian surface $M'$ of any other topological type would reveal that there are in fact at least $k$ distinct geodesics in $M'$ from $p$ to $q$ of length at most $kd$---a much tighter bound. A similar bound holds for geodesic loops at $p$. (See Theorem~1.4 from \cite{Nabutovsky_QuantitativeMorse}.)

This paper presents a significant improvement of the previously mentioned linear bounds for Riemannian 2-spheres as follows:
\begin{theorem}[Main Result]\label[theorem]{thm:TighterLinearBound}
    Let $p$ and $q$ be points on $\Sp$. Then for any metric $g$ on $\Sp$ and $k \geq 1$, the Riemannian 2-sphere $M = (\Sp,g)$ has at least $k$ distinct geodesics from $p$ to $q$ of length at most $8(k-1)d + \dist{p,q} \leq 8kd$,
    where $d = \diam(M)$. Moreover, for a set of metrics $g$ that is generic in the $C^m$ topology (for any $2 \leq m \leq \infty$), the tighter bound of $(7k - 6 - (k \bmod{2}))d + \dist{p,q} \leq 7kd$
    is satisfied.
    
    There are also at least $k$ distinct non-constant geodesic loops based at $p$ of length at most $6kd$. For a set of metrics $g$ that is generic in the $C^m$ topology for any $2 \leq m \leq \infty$, this can be further tightened to $(5k + (k \bmod{2}))d$.
\end{theorem}

(Note that if some periodic geodesic $\gamma$ passes through $p$, the geodesic loops based at $p$ from the statement of \cref{thm:TighterLinearBound} could be the iterates of $\gamma$, and may not be geometrically distinct. If $\gamma$ also passes through $q$, then the geodesics from $p$ to $q$ in the theorem statements could be iterates of $\gamma$ joined with with some arc of $\gamma$.)

The best known curvature-free bounds on the lengths of geodesics between given points will first be halved, by bounding the size of the ``gaps'' in the set of geodesic lengths. This will be established via a combination of the pigeonhole principle and the homology algebra structure of $\loops{pq}M$, the space of paths from $p$ to $q$. The bound will be lowered further using a new algorithm on the cut locus of $p$ in $M$. The algorithm constructs a ``sweep-out'' of $M$ by curves, whose lengths satisfy an upper bound which decreases when the metric of $M$ is ``more even.'' The rest of the proof applies the minimax principle, Lyusternik-Schnirelmann theory and Morse theory to $\loops{pq}M$, following prior literature. 

Our proof of \cref{thm:TighterLinearBound} will build on the arguments used in \cite{LinearBounds} and \cite{LinearBoundsSegments}, and will improve it by making its constructions more efficient, and by introducing some new ideas. Hence we will begin by summarizing the argument of \cite{LinearBoundsSegments} in \cref{sec:PriorWork}, and then outline our novel contributions in \cref{sec:Contributions}. The prior work is founded upon the general \emph{minimax principle} of Birkhoff and Hestenes \cite{BirkhoffHestenes_Minimax}, which manifests in our context as follows: for any Riemannian manifold $M$ and points $p,q \in M$, every homology class $h$ in the singular homology group $H_*(\loops{pq}M;\Q)$ is associated with some geodesic from $p$ to $q$. Moreover, its length is equal to the \emph{critical level} defined as
\begin{equation}
	\crit(h) = \inf_{Z \in h}\sup_{\alpha \in \supp(Z)} \length(\alpha),
\end{equation}
where $\supp(Z)$ denotes the support of the singular cycle $Z$, a set of paths that is the union of the images of each singular simplex in $Z$. Note that $\crit(h)$ depends on the length functional on $\loops{pq}M$, which in turn depends on the metric of $M$ and the choices of $p$ and $q$. Intuitively speaking, if $Z$ is a singular cycle representing a nontrivial homology class $h \in H_*(\loops{pq}M;\Q)$, then one may apply a curve-shortening process to all curves in $\supp(Z)$ at once, while fixing their endpoints. Since $h \neq 0$, the curves cannot all converge to the minimizing geodesic from $p$ to $q$. Instead, at least one of the curves must converge to a new geodesic from $p$ to $q$.\footnote{This flexible principle has been adapted in many other settings, involving other functionals and the homology of other spaces, to prove the existence of various ``minimal objects.'' These include the simple closed geodesics \cite{LyusternikSchnirelmann_ThreeGeodesics,Lyusternik_Book,LNR_PeriodicGeodesicsS2} and minimal hypersurfaces \cite{ColdingDeLellis_Minmax,MarquesNeves_MinimalVarieties}.}

Geodesics from $p$ to $q$ can be distinguished by two invariants: their lengths and, if the length functional on $\loops{pq}M$ is Morse (for instance, if $q$ is not conjugate to $p$ along any geodesic), their Morse indices. An approach common to both \cite{LinearBounds,LinearBoundsSegments} was to demonstrate that the homology classes of $\loops{pq}M$ ($\loops{pp}M$ in \cite{LinearBounds}, which we will henceforth denote by $\loops{p}M$) have critical levels that grow at most linearly, and to prove a dichotomy: either enough of them are distinct to satisfy the $22kd$ bound from \cite{LinearBoundsSegments}, or else many short geodesics from $p$ to $q$ can be constructed as local minima of the length functional. When the length functional is Morse, this dichotomy can be stated even more simply: \emph{either there will be many short geodesics of Morse index 0, or there will be many short geodesics with Morse indices 1, 2, 3 and so on.}

\subsection{Notation}

If $M$ is a Riemannian manifold, let $\diam(M)$ and $\inj(M)$ denote its diameter and injectivity radius respectively. For brevity, we may use the same symbol to refer to a path or its image; its precise meaning will be clear from the context. The concatenation of paths $\alpha$ followed by $\beta$ will be denoted by $\alpha \concat \beta$. The reversal of a path $\alpha$ will be written as $\rev{\alpha}$. The constant path at $p \in M$ will be denoted by $\const{p}$. If $\{\alpha_t\}_{t \in I}$ is an endpoint-fixing homotopy through curves of length at most $L$, then we may write it as $\alpha_t : \alpha_0 \homto{L} \alpha_1$. We may also write $\alpha_0 \homto{L} \alpha_1$ to simply indicate the existence of such a homotopy from $\alpha_0$ to $\alpha_1$.

The cardinality of a set $A$ will be denoted by $\card{A}$. Depending on the context, $\homot$ will indicate homotopy equivalence or diffeomorphism; $\iso$ will denote group isomorphism. We will also write $o(\epsilon)$ to indicate a value that can be made as small as desired by some choice of parameters in a relevant construction. Note that this differs from its usual designation as an asymptotically bounded quantity.

\subsection{Prior work}
\label{sec:PriorWork}

Fix two points $p$ and $q$ in $\Sp$. Let $M = (\Sp,g)$ be a Riemannian sphere. Rational homotopy theory may be used to show that $H_i(\loops{p}M;\Q) \iso \Q$ for all $i \geq 1$, and $H_i(\loops{p}M;\Q)$ is generated by the Pontryagin power $u^i$ of some generator $u \in H_1(\loops{p}M;\Q)$ \cite{Felix_RationalHomotopyTheory}. There is a homotopy equivalence $\varphi : \loops{p}M \to \loops{pq}M$ defined by concatenation with a fixed minimal geodesic from $p$ to $q$. This implies that $H_i(\loops{pq}M;\Q)$ is generated by homology classes $v^i = \varphi_*(u^i)$. Applying the minimax principle to the $v^0, v^1, v^2, \dotsc$ yields geodesics $\gamma_0, \gamma_1, \gamma_2, \dotsc$ where
\begin{equation}
    \label{eq:CriticalLevelsBoundU}
    \length(\gamma_i) = \crit(v^i) \leq \crit(u^i) + \dist{p,q} \leq \crit(u^i) + d.
\end{equation}

If the metric $g$ gives rise to a length functional that is Morse, then a standard argument shows that $\gamma_i$ has Morse index $i$, so the $\gamma_i$ must all be distinct.\footnote{For example, apply Theorem~17.3 from \cite{MorseTheory_Milnor} and cellular homology.} However, for a general metric $g$, it may be possible that $\gamma_i = \gamma_{i+1}$ and $\crit(v^i) = \crit(v^{i+1})$ for some values of $i$. On the other hand, Schwartz adapted the ideas of Lyusternik and Schnirelmann to prove that either $\crit(v^0) < \crit(v^2) < \crit(v^4) < \crit(v^6) < \dotsb$ and so on, or else if $\crit(v^{2m}) = \crit(v^{2m - 2})$, then there must be infinitely many geodesics from $p$ to $q$ of approximately that length. This is due to the cap product structure on $H_*(\loops{pq}M;\Q)$: there is a generator $c$ of $H^2(\loops{pq}M;\Q)$ such that $v^{2m} \frown c = v^{2m-2}$ for all $m \geq 1$. Lyusternik-Schnirelmann theory\footnote{An exposition of Lyusternik-Schnirelmann theory is available in Section 1 of \cite{Ballmann_LyusternikSchnirelmann}.} implies that if $a \frown c = a'$ for $0 \neq a \in H_*(\loops{pq}M;\Q)$, then either $\crit(a) > \crit(a')$ or else there will be an infinite sequence of geodesics from $p$ to $q$ whose lengths approach $\crit(a) = \crit(a')$. Unfortunately, this does not help to distinguish $\crit(v^{2m+1})$ from the other critical levels, because $v^{2m+1} \frown c = 0$ for all $m \geq 1$. Nevertheless, the cap product structure of the even powers of $v$ allows us to conclude that there are at least $k$ geodesic loops of length at most $\crit(v^{2k-2})$. Moreover, since the Pontryagin product concatenates loops and thus it adds lengths,
\begin{align}
    \crit(ab)
        &\leq \crit(a) + \crit(b) \text{ for all } a, b \in H_*(\loops{p}M;\Q) 
    \label{eq:PontryaginRelation_Loops} \\
    \implies \crit(v^{2k-2})
        &\leq (k-1)\crit(u^2) + \dist{p,q} 
    \label{eq:CriticalLevelsEvenBoundU2}\\
        &\leq (2k-2)\crit(u) + \dist{p,q}.  
    \label{eq:CriticalLevelsEvenBoundU}
\end{align}

In \cite{LinearBounds} and \cite{LinearBoundsSegments}, bounds on $\crit(u)$ and $\crit(u^2)$ were derived using \emph{sweep-outs}, defined as follows:
\begin{definition}[Sweep-out, induced meridians]
    A \emph{sweep-out} of $M$ is a continuous map $f : \Sp \to M$ of nonzero degree. The images of the meridians of $\Sp$ under $f$ are the \emph{induced meridians of $M$}.
\end{definition}

If $M$ admits a sweep-out with short induced meridians, we may describe it as a ``short sweepout''. The following result used by \cite{LinearBoundsSegments} would also imply bounds on $\crit(u)$ and $\crit(u^2)$. We state this result here in our language and sketch a short proof.

\begin{lemma}\label{lem:MeridiansBoundCriticalLevels}
    Let $M$ be a Riemannian 2-sphere with a chosen point $p \in M$. Suppose that $M$ has a sweep-out that maps a pole of $\Sp$ to $p$ so that the longest induced meridian of $M$ has length at most $L$ and the shortest induced meridian has length at most $l$. Then $\crit(u) \leq l + L$ and $\crit(u^2) \leq 2L$. 
\end{lemma}
\begin{proof}
    The induced meridians form a family of paths $\{\mu_t\}$, a shortest one being $\mu_{t_0}$. The lemma then follows from the observation that $u$ is represented by a 1-cycle consisting of loops $\mu_{t_0} \concat \rev{\mu_t}$, and that $u^2$ is represented by a 1-cycle consisting of loops $\mu_s \concat \rev{\mu_t}$.
\end{proof}

Hence, one might hope to bound $\crit(u)$ and $\crit(u^2)$ by finding some universal constant $K$ such that for any metric $g$ on $\Sp$, $(\Sp,g)$ admits a sweep-out that maps a pole to $p$ and has induced meridians of length at most $Kd$. However, a counterexample constructed by Y.~Liokumovich showed that such a universal constant $K$ does not exist: for any constant $K$, there exists an ``extremely rugged'' metric $g$ on $\Sp$ such that $(\Sp, g)$ has no sweep-out whose induced meridians are shorter than $Kd$ \cite{Liokumovich_LongSweepouts}. This counterexample was in turn based on an earlier construction, by S.~Frankel and M.~Katz, of an ``extremely rugged'' Riemannian 2-disk with a small diameter, and whose boundary is short but cannot be contracted to a point through short loops \cite{FrankelKatz}.

Nabutovsky and Rotman worked around this issue by establishing a dichotomy between the following scenarios for a Riemannian 2-sphere $M$:
\begin{itemize}
    \item $M$ has a short sweep-out, which gives short geodesics from $p$ to $q$ of length $\crit(v^0)$, $\crit(v^2)$, $\crit(v^4)$ and so on, where each $\crit(v^{2i})$ can be bounded in terms of the length of the longest induced meridian in the sweep-out.
    
    \item $M$ contains a \emph{convex annulus} or a simple periodic geodesic that contains $p$ and $q$. This convex annulus or periodic geodesic has infinitely many short geodesics winding around it.
\end{itemize}

(For simplicity, in the remainder of this Introduction, we will ignore the non-generic situation where a simple periodic geodesic passes through both $p$ and $q$.)

\begin{figure}[h]
    \centering
    \includegraphics{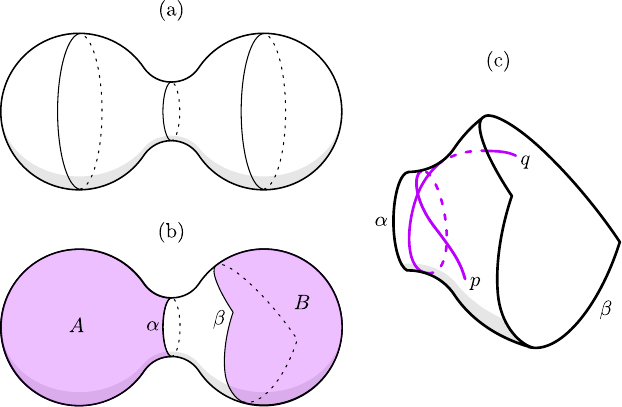}
    \caption{(a) A Riemannian 2-sphere $M$ in the shape of a dumbbell. (b) $\alpha$ is a periodic geodesic that bounds the shaded bottleneck disk $A$ (see \cref{def:BottleneckDisk}). $\beta$ is a loop composed of two geodesic arcs, and it also bounds a shaded bottleneck disk $B$. (c) The region between $\alpha$ and $\beta$ is a convex annulus. Any two points $p$ and $q$ in the annulus are endpoints of infinitely many geodesics that wind around the annulus, one of which is shown.}
    \label{fig:DumbbellBottleneckDisks}
\end{figure}

A convex annulus is a subspace $A$ of $M$ that is homeomorphic to a closed annulus such that $\partial A$ is \emph{convex to $A$} in the sense defined by C.\,B. Croke in Section 2 of \cite{Croke}. That is, $A$ must contain the minimal geodesics between all pairs of sufficiently close points in $\partial A$. (\Cref{fig:DumbbellBottleneckDisks}(c) illustrates an example of a convex annulus.) \cite{Croke}, \cite{LinearBounds} and \cite{LinearBoundsSegments} only consider the convexity of regions with piecewise geodesic boundary. To avoid unnecessary technicalities, this paper will also adopt the following more restrictive definition:

\begin{definition}[Convexity, concavity]
    \label[definition]{def:ConvexityConcavity}
    A region $A \subset M$ is \emph{convex} if $\partial A$ is piecewise geodesic and the internal angles at each vertex (measured from within $A$) are at most $\pi$. We will say that $A$ is \emph{convex except possibly at $x$} for some $x \in \partial A$ when $x$ is a vertex of $A$ and the angles at every other vertex are at most $\pi$. We will call $A$ \emph{concave} if its complement in $M$ is convex.
\end{definition}

\begin{figure}[h]
    \centering
    \includegraphics{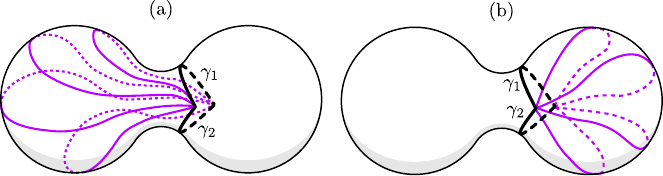}
    \caption{On the ``dumbbell'' Riemannian 2-sphere from \cref{fig:DumbbellBottleneckDisks}, $\gamma_1$ and $\gamma_2$ are minimizing geodesics with the same endpoints. The curves in each picture represent snapshots of a fixed-endpoint homotopy from $\gamma_1$ to $\gamma_2$ that must pass through long curves, due to the presence of the bottleneck disks in \cref{fig:DumbbellBottleneckDisks}(b).}
    \label{fig:DumbbellLongHomotopies}
\end{figure}

Observe that the complement of a convex annulus is the disjoint union of two concave open disks. The dichotomy between a short sweep-out and a convex annulus arises from the phenomenon that concave open disks can be shaped like ``bottlenecks'' that have short boundaries but large interiors. These ``bottlenecks'' obstruct short sweep-outs by forcing some curves that pass through concave open disks to be long (see \cref{fig:DumbbellLongHomotopies}). This inspires us to make the following definition:

\begin{definition}[Bottleneck disk]
    \label[definition]{def:BottleneckDisk}
    A \emph{bottleneck disk} is a concave subspace $D$ of a Riemannian 2-sphere $M$ that is homeomorphic to an open disk, such that $\length(\partial D) \leq 2\diam(M)$ and $\partial D$ is the union of at most two geodesics.\footnote{This condition is only present to reduce analytic technicalities when taking limits in later arguments.} (Occasionally we will allow $\partial D$ to be the union of at most four geodesics, but we will explicitly highlight these exceptions as they occur.)
\end{definition}

Using the language of bottleneck disks, we may formulate the dichotomy between short sweep-outs and convex annuli, as derived in 
\cite{LinearBounds} and \cite{LinearBoundsSegments}, in the following theorem.

\begin{theorem}[\cite{LinearBounds,LinearBoundsSegments}]
    \label[theorem]{thm:SweepoutOrConvexAnnulus}
    Let $M$ be a Riemannian 2-sphere containing points $p$ and $q$. Then either $M\setminus\{p,q\}$ contains 2 disjoint bottleneck disks, or $M$ admits a sweep-out that sends a pole to $p$ and induces meridians on $M$ of length at most $11d + o(\epsilon)$, where $d = \diam(M)$. If $p = q$, then the length bound can be tightened to $10d + o(\epsilon)$.
\end{theorem}

The proof of this theorem used the following techniques, which may be viewed as three variants of a discretized curve-shortening flow. We will also employ them later.
\begin{itemize}
    \item The \emph{Birkhoff curve shortening process for free loops} (BPFL) begins with any loop in a Riemannian manifold and produces a homotopy of loops, non-increasing in length, that converges either to a constant loop or a periodic geodesic. BPFL was first defined in \cite{Croke}.
    
    \item The \emph{Birkhoff curve shortening process for based loops} (BPBL) begins with a loop based at some point $x$ and produces a homotopy of loops based at $x$, non-increasing in length, that converges either to the constant loop at $x$ or a geodesic loop based at $x$. BPBL was first defined in \cite{LinearBounds}.
    
    \item The \emph{Birkhoff curve shortening process for segments} (BPS) begins with a path from $x$ to $y$ and produces a homotopy of paths from $x$ to $y$, non-increasing in length, that converges a geodesic from $x$ to $y$. BPS was first defined in \cite{LinearBoundsSegments}.
\end{itemize}

Detailed definitions of BPFL, BPBL, and BPS, as well as explanations for the convergence properties listed above, are available in Section~2 of \cite{LinearBoundsSegments}.\footnote{One advantage that BPBL and BPS have over BPFL in certain situations is that they allow one to control the location of the limiting geodesic loops or geodesic segments relative to $p$ and $q$. This property will help us.}

If $M$ contains a convex annulus $A$ that has short boundary components and that also contains both $p$ and $q$, then $A$ contains an infinite sequence of geodesics $\tau_1,\tau_2, \dotsc$ from $p$ to $q$ such that $\tau_k$ winds $k$ times around $A$ and has length at most $(2k + 4)d + o(\epsilon)$.\footnote{For details, see Section~2 of \cite{LinearBoundsSegments}, in particular Fig.~1(e).} These geodesic segments are obtained by applying BPS to short curves that represent different elements in the relative fundamental group $\pi_1(A,\{p,q\})$. Consequently, they are distinct. In a sense, the ``worst-case scenario'' occurs when $M$ does not contain such a convex annulus. In this scenario, Nabutovsky and Rotman used BPFL, BPBL and ``Gromov's pseudo-extension argument''\footnote{See Section~2, pages 415--417 of \cite{LinearBoundsSegments}.} in \cite{LinearBounds} and \cite{LinearBoundsSegments} to construct a sweep-out of $M$ whose longest and shortest induced meridians had lengths at most $11d + o(\epsilon)$ and $d + o(\epsilon)$ respectively. This, together with \cref{eq:CriticalLevelsEvenBoundU2}, the bound on $\crit(u^2)$ from \cref{lem:MeridiansBoundCriticalLevels}, and Lyusternik-Schnirelmann theory, implies that $M$ contains at least $k$ distinct geodesic segments from $p$ to $q$ of length at most $\crit(v^{2k-2}) \leq 22kd$. This yields the bound from \cite{LinearBoundsSegments}.

\subsection{Novel contributions}
\label{sec:Contributions}

\subsubsection*{Efficient constructions via a cut locus algorithm and radial convexity}

The $22kd$ and $20kd$ bound for geodesic segments and loops in \cite{LinearBoundsSegments} and \cite{LinearBounds} were derived using the same dichotomy between convex annuli and short sweep-outs. However, to require that the convex annulus contains two distinct points $p$ and $q$ is stricter than requiring it to contain just a single point $p = q$. Consequently, the construction of a short sweep-out of $M$ when $M\setminus\{p,q\}$ does not contain 2 disjoint bottleneck disks is more complicated in \cite{LinearBoundsSegments} ($p \neq q$), and yields a looser bound on the lengths of geodesic segments, as compared to the construction in \cite{LinearBounds} ($p = q$) and the bound it yielded for geodesic loops. In this paper we unify both cases by studying the relationship between the number of disjoint bottleneck disks contained in $M\setminus\{p\}$ and the lengths of induced meridians in sweep-outs. Thus we prove the following theorem, which strengthens \cref{thm:SweepoutOrConvexAnnulus} because if $M\setminus\{p\}$ contained 4 disjoint bottleneck disks, then $M\setminus\{p,q\}$ must contain 3 disjoint bottleneck disks.

\begin{theorem}
    \label[theorem]{thm:SweepoutOrBottleneckDisks}
    Let $M$ be a Riemannian 2-sphere and choose any $p \in M$. Then for $n \in \{2,4\}$, either $M \setminus \{p\}$ contains $n$ disjoint bottleneck disks (whose boundaries may be the unions of at most four geodesics), or $M$ admits a sweepout that maps one pole to $p$ such that the longest induced meridian on $M$ has length at most $(n+3)d + o(\epsilon)$ and the shortest induced meridian has length at most $d + o(\epsilon)$.
\end{theorem}

The $n = 2$ case is used to prove the bound for geodesic loops in our main result, while the $n = 4$ case is used to prove the bound on geodesics from $p$ to $q$.\footnote{Our proof of \cref{thm:SweepoutOrBottleneckDisks} can be extended to prove the theorem for all even integers $n \geq 2$, which may be of independent interest.} The $5d + o(\epsilon)$ bound for the $n = 2$ case in \cref{thm:SweepoutOrBottleneckDisks} is tighter than the $7d + o(\epsilon)$ bound for the $n = 4$ case, and this causes our bound for geodesic loops in our main result to be tighter than the bound for geodesics from $p$ to $q$.

The improvements of the bounds from \cref{thm:SweepoutOrBottleneckDisks} over those in \cref{thm:SweepoutOrConvexAnnulus} are partly due to more efficient constructions that exploit the cut locus of $p$ in $M$, which is a tree when $M$ has an analytic metric \cite{Myers_AnalyticCutLocus}. We will present a recursive algorithm on the cut locus that ``searches'' down the tree for bottleneck disks; if it cannot find enough of them, the algorithm can construct fixed-endpoint homotopies between certain pairs of minimizing geodesics that pass through short curves. These homotopies are then are patched together into short sweep-outs. We note that a similar idea is present in \cite{Beach_FreeLoopSpace}.

The fixed-endpoint homotopies in the cut locus algorithm are constructed efficiently with the help of the notion of \emph{radially convex sets},\footnote{In the context of Euclidean spaces, radially convex sets are also commonly known as ``star-shaped sets''.} which we generalize to Riemannian manifolds as follows:

\begin{definition}[Radial convexity]
    \label[definition]{def:RadialConvexity}
    Given a subspace $A$ of a Riemannian manifold $M$ and some point $x \in A$, we say that $A$ is \emph{radially convex from $x$} if it contains every minimizing geodesic from $x$ to a point in its interior.
\end{definition}

We show that in the key lemmas that were used in most of the constructions from \cite{LinearBounds} and \cite{LinearBoundsSegments}, many of the subspaces $A$ of $M$ considered in the lemmas are in fact radially convex from $p$, and that applying BPBL to $\partial A$ yields curves that bound regions that remain radially convex from $p$. This ensures that the property of radial convexity is compatible with the use of BPBL in \cite{LinearBounds} and \cite{LinearBoundsSegments}, allowing us to derive tighter bounds using similar constructions.

\subsubsection*{Avoiding bound-doubling using subadditivity}

In \cite{LinearBounds} and \cite{LinearBoundsSegments}, the possibility that some homology classes in $v^0, v^1, v^2, \dotsc$ might give rise to the same geodesic was handled using Lyusternik-Schnirelmann theory. However, that caused the $22kd$ bound from \cite{LinearBoundsSegments} to be twice as large as those in \cref{thm:SweepoutOrConvexAnnulus}. In this paper we avoid doubling the bounds derived from \cref{thm:SweepoutOrBottleneckDisks} by using more of the relations from \cref{eq:PontryaginRelation_Loops} than were used to derive \cref{eq:CriticalLevelsEvenBoundU2,eq:CriticalLevelsEvenBoundU}. The motivation for our approach is apparent from the simpler case where $p = q$ (and thus $u^i = v^i$ for all $i$): suppose that we can prove that $\crit(u) \leq L$, so that by  \cref{eq:CriticalLevelsEvenBoundU}, $\length(\gamma_{2k}) \leq 2k\crit(u) \leq 2kL$. In the ``worst case,'' the geodesic loops $\gamma_0, \gamma_2, \gamma_4, \dotsc$ might have lengths exactly 0, $2L$, $4L$, and so on. The concern behind the bound-doubling was that perhaps $\gamma_{2k} = \gamma_{2k+1}$ for some or even all values of $k$. However, this cannot happen for any $k \geq 0$, otherwise \cref{eq:PontryaginRelation_Loops} would imply a contradiction: $(2k+2)L = \length(\gamma_{2k+2}) = \crit(u^{2k+2}) \leq \crit(u^{2k+1}) + \crit(u) = \crit(u^{2k}) + \crit(u) \leq (2k + 1)L$. Intuitively, if $\gamma_{2k} = \gamma_{2k+1}$, then $\gamma_{2k+2}, \gamma_{2k+3}, \gamma_{2k+4}, \dotsc$ must satisfy even tighter length bounds to prevent gaps larger than $L$ from appearing in the set $\{\crit(u^0), \crit(u^1), \crit(u^2), \dotsc\}$.
In the general case, we will use the Pigeonhole Principle to prove that the set of values $\{\crit(v^0), \crit(v^1), \crit(v^2), \dotsc\}$ cannot have ``gaps'' larger than $\crit(u)$. This will allow us to prove the following:

\begin{proposition}\label[proposition]{prop:LinearBoundFirstCriticalLevel}
    For any Riemannian sphere $M$ and points $p, q \in M$, let $u$ generate $H_1(\loops{p}M;\Q)$. Then there are at least $k$ distinct geodesics from $p$ to $q$ of length at most $(k-1)\crit(u) + \dist{p,q} \leq (k-1)\crit(u) + d$. There are also at least $k$ distinct non-constant geodesic loops based at $p$ of length at most $k\crit(u)$.
\end{proposition}

The above results can be collected into a proof of \cref{thm:TighterLinearBound} in a manner outlined as follows. When $M\setminus\{p\}$ does not contain 4 disjoint bottleneck disks, the $8(k-1)d + \dist{p,q}$ bound on lengths of geodesic segments from $p$ to $q$ in \cref{thm:TighterLinearBound} follows from the $n = 4$ case of \cref{thm:SweepoutOrBottleneckDisks}, \cref{lem:MeridiansBoundCriticalLevels}, as well as \cref{prop:LinearBoundFirstCriticalLevel}. The tighter bound for a generic set of metrics follows from the fact that for such a set of metrics, the length functional is Morse, so the geodesics $\gamma_0, \gamma_1, \gamma_2, \dotsc$ would have Morse indices 0, 1, 2, and so on, implying that they are distinct. When $M\setminus\{p\}$ contains 4 disjoint bottleneck disks, then there would have to be a convex annulus in $M$ that contains both $p$ and $q$. This implies that a much tighter bound holds as explained in \cref{sec:PriorWork}. A similar chain of derivations yields the bounds on the lengths of geodesic loops based at $p$.

\subsection{Organizational structure}

In \Cref{sec:SweepoutOrBottleneckDisks} we will present an algorithm that operates on a cut locus tree to construct short fixed-endpoint homotopies in a Riemannian 2-sphere $M$, when $M$ does not contain many disjoint bottleneck disks. A certain proposition (\cref{prop:MinBigon1ConcaveDiskOrHomotopy}) will be assumed and applied repeatedly in the algorithm, but its more technical proof will be deferred to \cref{sec:RadialConvexity}. The proposition helps to bound the lengths of the curves in the homotopies. Assuming this proposition, we will also prove \cref{thm:SweepoutOrBottleneckDisks}.

In \cref{sec:RadialConvexity} we will prove that radial convexity is preserved by BPBL. We will then adapt certain constructions from \cite{LinearBounds} and \cite{LinearBoundsSegments}, while exploiting radial convexity, to prove \cref{prop:MinBigon1ConcaveDiskOrHomotopy}. In \cref{sec:Combinatorics} we will use a generalization of the relation in \cref{eq:PontryaginRelation_Loops} and the pigeonhole principle to prove \cref{prop:LinearBoundFirstCriticalLevel}, allowing us to avoid doubling the bound like in the prior literature. In \cref{sec:MorseHomology} we will derive tighter length bounds in Riemannian 2-spheres with Morse length functionals, and finally tie all of our results together prove \cref{thm:TighterLinearBound}.

\subsection*{Acknowledgements}

The author would like to express his gratitude to his academic advisors, Regina Rotman and Alexander Nabutovsky, for suggesting this research topic and for valuable discussions. The author would also like to thank Angela Wu, whose unpublished notes motivated this work.

\section{Short sweep-outs or many pairwise disjoint bottleneck disks}
\label{sec:SweepoutOrBottleneckDisks}

To construct sweep-outs on a Riemannian 2-sphere $M$, it helps to decompose it into simpler pieces and try to ``sweep out'' those pieces, or look for bottleneck disks in each piece that might obstruct short sweep-outs. If $M$ has a generic analytic metric, then the cut locus of a point $p \in M$ is a tree \cite{Buchner_CutLocus,Myers_AnalyticCutLocus}. The minimizing geodesics from $p$ to a given vertex of the tree divide $M$ into \emph{minimizing geodesic bigons}, defined as follows. These bigons form a convenient decomposition of $M$ for the construction of sweep-outs.

\begin{definition}[Minimizing geodesic bigon]
    Given a Riemannian 2-sphere, a \emph{minimizing geodesic bigon} is a subspace homeomorphic to a closed disk whose boundary is the union of two minimizing geodesics that touch only at their endpoints, $x$ and $y$. The points $x$ and $y$ are called the \emph{vertices} of the minimizing geodesic bigon, and are uniquely defined if and only if the boundary of the bigon is not a periodic geodesic.
\end{definition}

The following result shows that if a minimizing geodesic bigon---a ``piece'' of $M$---does not contain a bottleneck disk, then this ``piece'' can be ``swept out'' by short curves.

\begin{proposition}
    \label[proposition]{prop:MinBigon1ConcaveDiskOrHomotopy}
    Let $B$ be a minimizing geodesic bigon in a Riemannian sphere $M$ whose boundary is the union of two minimizing geodesics $\gamma_1$ and $\gamma_2$. Then either $B$ contains a bottleneck disk or there is a homotopy $\gamma_1 \homto{5\diam(M)} \gamma_2$ through curves in $B$.
\end{proposition}

Lemmas~3.2, 3.4 and 3.5 from \cite{LinearBoundsSegments} imply a statement almost identical to \cref{prop:MinBigon1ConcaveDiskOrHomotopy}, except with the looser length bound of $6\diam(M)$. In \cref{sec:RadialConvexity} we will use the concept of radial convexity to prove the tighter bound in \cref{prop:MinBigon1ConcaveDiskOrHomotopy} However, a series of technical constructions are required for the proof, so we defer it to \cref{sec:RadialConvexity} for the sake of first explaining how the proposition can be applied to prove \cref{thm:SweepoutOrBottleneckDisks}.

\begin{proposition}
    \label[proposition]{prop:MinBigon2ConcaveDisksOrHomotopy}
    Let $g$ be a generic analytic metric on $\Sp$, and let $B \subset M = (\Sp,g)$ be a minimizing geodesic bigon whose boundary is the union of two minimizing geodesics $\gamma^-$ and $\gamma^+$. Then either $B$ contains 2 disjoint bottleneck disks or there is a homotopy $\gamma^- \homto{7\diam(M)} \gamma^+$ through curves in $B$.
\end{proposition}
\begin{proof}
    Let $x = \gamma^+(0)$ and $y = \gamma^+(1)$ be the vertices of $B$, and let $d = \diam(M)$. Using a recursive algorithm, we will attempt to contruct a contracting homotopy \begin{equation}
        \label{eq:ContractBigonBoundary}
        \gamma^- \concat \rev{\gamma^+} \homto{6d} \const{x},
    \end{equation}
    which we can then extend to a homotopy $\gamma^- \homto{3d} \gamma^- \concat \rev{\gamma^+} \concat \gamma^+ \homto{7d} \const{x}  \concat \gamma^+ \homto{d} \gamma^+$, where the first homotopy in the sequence passes through curves that extend $\gamma^-$ from $y$ along part of $\im\gamma^+$, and then backtrack to $y$. The recursive algorithm will fail only if we find two disjoint bottleneck disks in $B$. 
    
    The algorithm can be described as follows. Since $g$ is a generic analytic metric, the cut locus $C(x)$ of $x$ in $M$ is a tree with a vertex at $y$ \cite{Myers_AnalyticCutLocus,Buchner_CutLocus}. Moreover, we can assume that except possibly for $y$, the vertices of the tree are \emph{leaves} with degree 1 or \emph{interior vertices} with degree 3. If $C(x) \cap B = \{y\}$ then an analysis of the exponential map at $x$ reveals a homotopy $\gamma^- \homto{d} \gamma^+$ through minimizing geodesics. This can be extended to a contracting homotopy $\gamma^- \concat \rev{\gamma^+} \homto{2d} \gamma^+ \concat \rev{\gamma^+} \homto{2d} \const{x}$  as required in \cref{eq:ContractBigonBoundary}, where the second homotopy contracts loops along $\im\gamma^+$. Therefore we may proceed by assuming that $C(x) \cap B$ is a binary tree rooted at $y$, where $y$ can have either one or two child vertices in the tree. We will deal with the case where $y$ has only one child $z$, as the other case can be handled analagously. There is a continuous family of minimizing geodesic bigons $B_t \subset B$, where $t \in [0,1)$, with one vertex at $x$ and the other vertex on the edge $yz$ of the cut locus. Let $\partial B_t = \gamma_t^- \cup \gamma_t^+$, where each $\gamma_t^\pm$ is a minimizing geodesic starting from $x$ which lies on the same side of $yz$ as $\gamma^\pm$. (See \cref{fig:CutLocusAlgorithm}(a) for an illustration of the notation.) Let $\gamma_t^\pm$ tend to a minimizing geodesic $\gamma_1^\pm$, giving us a homotopy
    \begin{equation}
        \label{eq:HomotopyAlongEdge}
        \gamma^- \concat \rev{\gamma^+} \homto{2d} \gamma_1^- \concat \rev{\gamma_1^+}
    \end{equation}
    If $\gamma_1^- = \gamma_1^+$ then we may finish with a homotopy $\gamma_1^- \concat \rev{\gamma_1^+} \homto{2d} \const{x}$ that contracts the loop along $\im\gamma_1^+$. Otherwise, $\gamma_1^- \cup \gamma_1^+$ bounds a minimizing geodesic bigon $B_1 = \bigcap_tB_t$. Note that $z$ is either a leaf of the tree or else it must have exactly two child vertices. 
    \begin{itemize}
        \item Suppose that $z$ happens to be a leaf of the tree. Then there will be a homotopy $\gamma_1^- \concat \rev{\gamma_1^+} \homto{2d} \const{x}$ similar to the situation where $C(x) \cap B = \{y\}$. Combining this homotopy with the one from \cref{eq:HomotopyAlongEdge} gives the required contracting homotopy in \cref{eq:ContractBigonBoundary}.
        
        \item Suppose that $z$ is an interior vertex with two child vertices. Then there must also be a third minimizing geodesic $\gamma_1^0$ from $x$ to $z$ that is different from $\gamma_1^\pm$. Let $B_1^\pm$ be the minimizing geodesic bigon bounded by $\gamma_1^\pm \cup \gamma_1^0$. If both $B_1^-$ and $B_1^+$ contain a bottleneck disk, then we have satisfied the conclusion of the proposition and we can terminate the algorithm. Otherwise, we may assume without loss of generality that $B_1^-$ does not contain a bottleneck disk, so \cref{prop:MinBigon1ConcaveDiskOrHomotopy} gives a homotopy $\gamma_1^- \homto{5d} \gamma_1^0$, which can be extended to a homotopy $\gamma_1^- \concat \rev{\gamma_1^+} \homto{6d} \gamma^0 \concat \rev{\gamma^+}$. Combining this with the homotopy from \cref{eq:HomotopyAlongEdge} gives a homotopy $\gamma^- \concat \rev{\gamma^+} \homto{6d} \gamma^0 \concat \rev{\gamma^+}$. Now we can return to the start of the algorithm and apply it to $B_1^+$ (see \cref{fig:CutLocusAlgorithm}(b)).
    \end{itemize}

    \begin{figure}[h]
        \centering
        \includegraphics{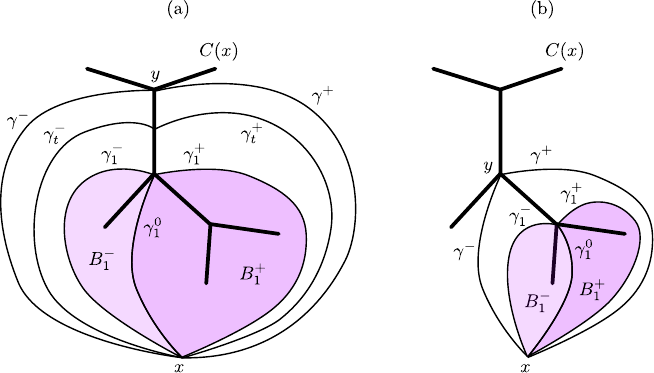}
        \caption{(a) An illustration of what the first iteration of the cut locus algorithm could look like. (b) An illustration of what the second iteration could look like, when applied to the bigon labeled $B_1^+$ in (a), but with a new label for each element that describes its role in this new iteration. For example, $\gamma_1^0$ from (a) now plays the role of $\gamma^-$.}
    \label{fig:CutLocusAlgorithm}
    \end{figure}
    
    This recursion will continue until we have found a chain of homotopies that combine into the the one required by \cref{eq:ContractBigonBoundary}, or else we will find two disjoint bottleneck disks in $B$ at some point.
\end{proof}

Now we are ready to prove \cref{thm:SweepoutOrBottleneckDisks}, as long as we assume \cref{prop:MinBigon1ConcaveDiskOrHomotopy}.

\begin{proof}[Proof of \cref{thm:SweepoutOrBottleneckDisks}]
    First we prove the theorem in the situation where $M$ has a generic analytic metric---and where the bottleneck disks are required to have boundaries consisting of at most two geodesics as usual---and then extend the result to all smooth metrics, where the bottleneck disks are allowed to have boundaries consisting of at most four geodesics.
    
    Assume that $M$ has a generic analytic metric. For the $n = 2$ case, suppose that $M \setminus \{p\}$ does not contain 2 disjoint bottleneck disks. The required sweep-out can be constructed following the procedure from Section~2 of \cite{LinearBoundsSegments}, with a single modification as follows. The crucial step is to find fixed-endpoint homotopies between any pair of distinct minimizing geodesics from $p$ to some point $x$ on its cut locus. These two minimizing geodesics divide $M$ into two minimizing geodesic bigons $B$ and $B'$. But either $B$ or $B'$ must contain no bottleneck disk, in which case \cref{prop:MinBigon1ConcaveDiskOrHomotopy} supplies the required homotopy through curves of length at most $5d$. An analogous argument works for the $n = 4$ case, except that we assume that $M \setminus \{p\}$ does not contain 4 disjoint bottleneck disks, so either $B$ or $B'$ must contain at most one bottleneck disk, and thus \cref{prop:MinBigon2ConcaveDisksOrHomotopy} supplies the required homotopies through curves of length at most $7d$.
    
    To extend the result to all smooth metrics on $\Sp$, we adapt some of the arguments used at the end of the proof of Theorem~0.1 in \cite{LinearBoundsSegments}. Let $g$ be an arbitrary smooth metric on $\Sp$, and find a sequence $g_1, g_2, \dotsc$ of generic real analytic metrics on $\Sp$ that converges to $g$ in the $C^2$ topology. Write $M_i = (\Sp,g_i)$ and $M = (\Sp,g)$. The sequence of diameters $d_i = \diam(M_i)$ converges to $d = \diam(M)$. By passing to a subsequence, we may assume that we are in either of the following situations for $n \in \{2,4\}$, because we have proven the theorem for generic analytic metrics:
    \begin{itemize}
        \item Each $M_i$ admits a sweepout that sends one pole to $p$ and has longest and shortest induced meridians of bounded length. As outlined in the proof of Theorem~0.1 in \cite{LinearBoundsSegments}, one of those sweepouts can be converted to a sweep-out of $M$ satisfying the same bounds, up to an $o(\epsilon)$ error.
        
        \item Every $M_i\setminus\{p\}$ contains $n$ disjoint bottleneck disks $D_{i1}, \dotsc, D_{in}$. We will use this to show that $M\setminus\{p\}$ also contains $n$ disjoint bottleneck disks whose boundaries may be the unions of at most four geodesics. Parametrize each $\partial D_{ij}$ by a curve $\alpha_{ij}$. The Arzel\`{a}-Ascoli theorem allows us to pass to a subsequence of $\{M_i\}$ such that for each $1 \leq j \leq n$, the sequence $\alpha_{1j}, \alpha_{2j}, \dotsc$ converges to a closed curve $\alpha_{\infty j}$ in the compact-open topology. The limit curves $\alpha_{\infty 1}, \dotsc, \alpha_{\infty n}$ must have length at most $2d$, and each must be the union of at most two geodesics. We will modify the limit curves until they are simple and bound pairwise disjoint bottleneck disks that do not contain $p$ and that are bounded by at most four geodesics.
        
        The possibility that $\alpha_{\infty j}$ could be the constant curve or the same geodesic travelled twice in opposite directions can be ruled out by applying Cheeger's inequality\footnote{See Corollary~2.2 in \cite{Cheeger_InjectivityRadius}.} in a way similar to the proof of Theorem 0.1 in \cite{LinearBoundsSegments}, and by considering that each $D_{ij}$ is concave and cannot lie in a disk with radius less than $\inj(M)$.
    
        $\alpha_{\infty j}$ cannot ``cross'' itself because it is the limit of simple curves, but it may ``touch'' itself at its vertices without crossing itself, as depicted in \cref{fig:SelfTouching}(a). $\alpha_{\infty j}$ can be modified into a simple closed curve by ``chipping off a corner'' at every ``self-touching'' vertex $v$ as shown in \cref{fig:SelfTouching}(b). This modification can be sufficiently small to avoid introducing new intersections between $\alpha_{\infty j}$ and the other curves $\alpha_{\infty k}$; otherwise, $v$ would also have to be a vertex of $\alpha_{\infty k}$, which would contradict the concave nature of one of the disks $D_{ik}$. The modification can also be made so that the resulting curves bound pairwise disjoint bottleneck disks $D_{\infty 1}, \dotsc, D_{\infty n}$ that do not contain $p$, such that each $\partial D_{\infty j}$ is the union of at most four geodesics (see \cref{fig:SelfTouching}(b)).
    \end{itemize}
    
    \begin{figure}
        \centering
        \includegraphics{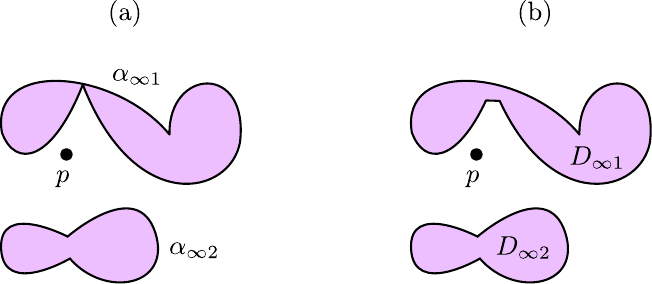}
        \caption{An illustration of possible limit curves $\alpha_{\infty j}$ where $n = 2$. The self-touching $\alpha_{\infty1}$ in (a) can be modified by ``chipping off a corner'' as shown in (b). The shaded regions in (b) represent bottleneck disks.}
        \label{fig:SelfTouching}
    \end{figure}
\end{proof}

\section{Using Radial Convexity to Prove Proposition~\ref{prop:MinBigon1ConcaveDiskOrHomotopy}}
\label{sec:RadialConvexity}

For this section, fix a Riemannian 2-sphere $M$ and a point $p \in M$. Let $d = \diam(M)$. We will prove \cref{prop:MinBigon1ConcaveDiskOrHomotopy} by proving more efficient versions of Lemmas~3.2, 3.4 and 3.5 from \cite{LinearBoundsSegments} that exploit radial convexity (\cref{def:RadialConvexity}). To modify these lemmas and their proofs in a self-contained fashion that uses our language of bottleneck disks, we will present our own proofs of some of them in this section.

\subsection{Properties of radial convexity}

Let us begin by establishing some useful properties of radial convexity. The relevance of radial convexity in constructions involving minimizing geodesics is apparent from the following lemma.

\begin{figure}[h]
    \centering
    \includegraphics{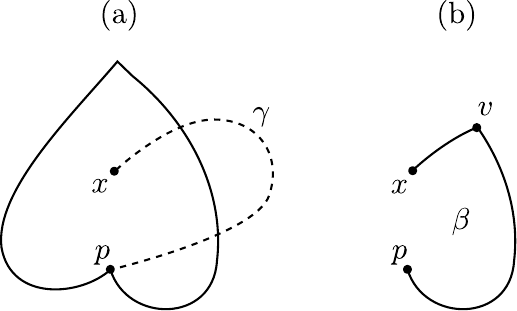}
    \caption{An illustration of the proof by contradiction for the radial convexity of minimizing bigons in the proof of \cref{lem:MinBigonRadiallyConvex}.}
    \label{fig:MinBigonRadiallyConvex}
\end{figure}

\begin{lemma}\label{lem:MinBigonRadiallyConvex}
    Each minimizing geodesic bigon is radially convex from both vertices.
\end{lemma}
\begin{proof}
    Assume for the sake of contradiction that $B$ is a minimizing geodesic bigon that is not radially convex from one of its vertices $p$. Thus the interior of $B$ contains some $x$ such that some minimizing geodesic $\gamma$ from $p$ does not stay in $B$. (See \cref{fig:MinBigonRadiallyConvex}(a).) Let $\gamma$ have length $\ell$. We will modify $\gamma$ into a piecewise geodesic path $\beta$ from $p$ to $x$ of length $\ell$ that stays in $B$, and then shorten it by ``cutting'' one of its corners while stay in $B$. This would contradict the fact that there can be no path from $p$ to $x$ that is shorter than $\ell$.
    
    Parametrize $\partial B$ by a constant-speed curve $\alpha : [0,1] \to \partial B$ that starts and ends at $p$. As a result, $\alpha(1/2) = x$. Since $\gamma$ must exit $B$ at some point, there must exist some $s, t > 0$ so that $\gamma(s) = \alpha(t)$; choose the greatest value of $s$ for this, so that $\gamma|_{[s,1]} \subset B$. If $t \leq 1/2$ then the fact that $\gamma$ is a minimizing geodesic implies that $\length(\gamma|_{[0,s]}) \leq \length(\alpha|_{[0,t]})$. The inequality in the other direction is a consequence of the fact that $\alpha|_{[0,t]}$ is a minimizing geodesic. Hence we have a path $\beta = \alpha|_{[0,t]} \concat \gamma|_{[s,1]} \subset B$ of length $\ell$. If $t > 1/2$ then a similar argument involving the minimizing geodesic $\alpha|_{[1/2,1]}$ gives the same result.
    
    Now we have a piecewise geodesic path $\beta = \alpha' \concat \gamma' \subset D$ of length $\ell$ from $p$ to $x$, where $\alpha' \subset \alpha$ and $\gamma' \subset \gamma$. Moreover, $\beta$ is defined in a way so that it must spend nonzero time on $\alpha$ and $\gamma$, and the point at which $\beta$ switches from $\alpha$ to $\gamma$ must be a vertex $v$ which has angle less than $\pi$ on one side. (See \cref{fig:MinBigonRadiallyConvex}(b).) Since there can be no path from $p$ to $x$ that is shorter than $\ell$, this path $\beta$ cannot be shortened within $B$ by ``cutting'' a corner near $v$, which forces $v = x$. That is, $\alpha'$ has to be either $\alpha|_{[0,1/2]}$ or $\alpha|_{[1/2,1]}$. However, now either $\alpha|_{[0,1/2]} \concat \gamma'$ or $\alpha|_{[1/2,1]} \concat \gamma'$ can be shortened by cutting a corner inside $B$, contradicting the fact that there can be no path from $p$ to $x$ that is shorter than $\ell$.
\end{proof}

The following lemma follows immediately from the definition of radial convexity.

\begin{lemma}\label{Lem:RadiallyConvexIntersection}
    Let $\{A_j\}_{j \in J}$ be a family of sets containing some $x \in M$ that are all radially convex from $x$. Then $\bigcap_jA_j$ is also radially convex from $x$.
\end{lemma}

The following lemma states some properties of BPBL that follow from arguments in \cite{LinearBoundsSegments} and \cite{Croke}.

\begin{lemma}
    \label{lem:BPBLProperties}
    Let $D \subset M$ be a closed disk that is convex except possibly at some point $x \in \partial D$. Suppose that BPBL applied to $\partial D$ while fixing $x$ yields the family of curves $\{\alpha_t\}$. Then all of these curves lie in $D$. Furthermore, let $t_0 \leq \infty$ be the infimum of all $t$ such that $\alpha_t$ self-intersects. The closed disks $\{D_t\}_{t < t_0}$ bounded by the simple loops $\{\alpha_t\}_{t < t_0}$ are all convex except possibly at $x$, and form a descending filtration: $D_t \subset D_s$ for all $s < t < t_0$. Finally, if $t_0 < \infty$ then $\bigcap_{t < t_0}D_t$ is the union of two convex closed disks that intersect only at $x$. The boundary of the union would then be $\alpha_{t_0}$.
\end{lemma}
\begin{proof}
    The claim that $D_t$ is convex except possibly at $x$ whenever $t < t_0$ follows from adapting Croke's proof of Lemma~2.2 in \cite{Croke}. The rest of the lemma follows from arguments of Nabutovsky and Rotman in their proof of Lemma~3.5 in \cite{LinearBoundsSegments}.
\end{proof}

Our next result demonstrates that under some conditions, radial convexity is preserved by BPBL.

\begin{figure}[h]
    \centering
    \includegraphics{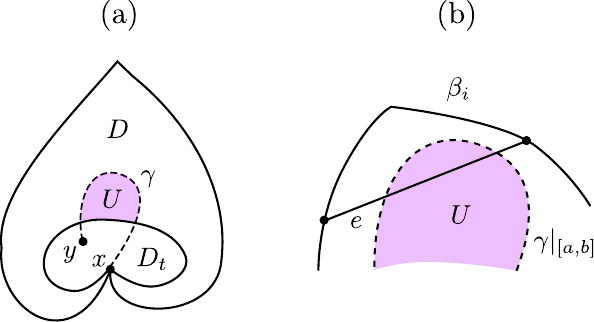}
    \caption{An illustration of how in \cref{lem:RadiallyConvexBCSP}, a minimizing geodesic whose endpoints are in $D_t$ but does not stay in $D_t$ would ``block'' BPBL. See the details in the proof.}
    \label{fig:RadiallyConvexBCSP}
\end{figure}

\begin{lemma}\label{lem:RadiallyConvexBCSP}
    Let $D$ be a closed disk that is convex except possibly at some point $x \in \partial D$, and is radially convex from $x$. Suppose that BPBL applied to $\partial D$ while fixing $x$ yields the curves $\alpha_t$, which are simple for all $t < t_0$ for some $t_0$. Let $D_t \subset D$ be a closed disk bounded by $\alpha_t$. Then $D_t$ is radially convex from $x$ for all $t < t_0$, and so is $\bigcap_{t < t_0}D_t$.
\end{lemma}
\begin{proof}
    It suffices to show that $D_t$ is radially convex from $x$ for all $t < t_0$, since the remaining claim will then follow from \cref{Lem:RadiallyConvexIntersection}. 
    Since $D$ is convex except possibly at $x$, \cref{lem:BPBLProperties} guarantees that for all $s < t$ we have $D_t \subset D_s \subset D$. Suppose for the sake of contradiction that $D_t$ is not radially convex from $x$ for some $t < t_0$. Then the interior of $D_t$ must contain some $y$ such that some minimizing geodesic $\gamma$ from $x$ to $y$ does not stay in $D_t$. However, $y$ must lie in the interior of the radially convex $D$, so $\gamma \subset D$. Let $\gamma|_{[a,b]}$ be a segment of $\gamma$ that lies outside $D_t$ except for its endpoints $\gamma(a)$ and $\gamma(b)$. We will show that such a segment would obstruct BPBL from reaching $\alpha_t$, yielding a contradiction.
    
    Let $U$ be the open disk bounded by $\gamma|_{[a,b]}$ and $\partial D_t$. (See \cref{fig:RadiallyConvexBCSP}(a).) Recall that BPBL on $\alpha$ fixing $p$ proceeds by subdividing $
    \alpha = \beta_0$ into segments whose lengths we can choose to be less than $\inj(M)/2$, joining the midpoints of those segments via geodesic arcs into a ``midpoint polygon'' $\beta_1$, interpolating between $\beta_0$ and $\beta_1$ via a length non-increasing homotopy, and then recursively applying this process to $\beta_1$. The result is a sequence of homotopies interpolating between the curves $\beta_0, \beta_1, \beta_2, \dotsc$ and so on.\footnote{For a detailed definition, see Section~2 in \cite{LinearBoundsSegments}.} Choose the smallest index $i$ such that $\beta_i$ is disjoint from $U$ but the interpolating homotopy from $\beta_i$ to $\beta_{i+1}$ contains a curve $\alpha_s$ that intersects $U$. We may choose $s$ to be less than $t$. That is, $\gamma$ must intersect the interior of a small geodesic triangle whose edges are minimizing geodesics, two of whom are part of $\beta_i$ while the last edge $e$ is part of $\alpha_s$. (See \cref{fig:RadiallyConvexBCSP}(b).) This means that $\gamma$ must intersect only $e$, and since the $e$ must be shorter than the injectivity radius, it must be part of the minimizing geodesic $\gamma$. This contradicts the assumption that $\alpha_s$ intersects $U$.
\end{proof}

\subsection{Geometric constructions that yield bottleneck disks or homotopies through short curves}

\Cref{lem:RadiallyConvexToPoint,lem:MinBigonBCSPSimple,lem:MinBigonBCSPSelfCross} below encapsulate the main geometric constructions that quantify the degree to which bottleneck disks obstruct homotopies through short curves---that is, how short the curves in a homotopy could be in the absence of bottleneck disks. Together they will help us prove \cref{prop:MinBigon1ConcaveDiskOrHomotopy}.

The following lemma can be derived using the proof of Lemma~3.4 from \cite{LinearBoundsSegments}, but with a better bound as explained in the subsequent remark.

\begin{lemma}
    \label{lem:RadiallyConvexToPoint}
    Let $D \subset M$ be a closed disk that contains $p$, is radially convex from $p$, and is convex. Suppose that applying BPFL to $\partial D$ makes it converge to a point $y \in D$. Then there is a contraction of $\partial D$ that passes through curves in $D$ and based at $p$ of length at most $L + 2d$, where $L = \length(\partial D)$.
\end{lemma}

\begin{remark}
    The length bound of $L + 2d$ in \cref{lem:RadiallyConvexToPoint} is tighter than the $2L + 2d$ bound in \cite{LinearBoundsSegments} because the radial convexity of $D$ from $p$ guarantees that $\dist[D]{y,p} \leq d$, where $\operatorname{dist}_D$ denotes the intrinsic distance function of $D$. The proof of Lemma~3.4 from \cite{LinearBoundsSegments} relied on the looser bound of $\dist[D]{y,p} \leq d + L/2$.
\end{remark}

\begin{lemma}\label{lem:MinBigonBCSPSimple}
    Let $B$ be a minimizing geodesic bigon with one vertex at $p$, and with an angle of at most $\pi$ at the other vertex when measured from the inside. Parametrize $\partial B$ as a simple loop $\alpha : [0,1] \to \partial B$ based at $p$. Suppose that the BPBL on $\alpha$ that fixes $p$ only passes through simple curves. Then either $B$ contains a bottleneck disk or $\alpha$ can be contracted to $p$ through curves in $B$ that are based at $p$ and have length at most $\length(\alpha) + 2d$.
\end{lemma}
\begin{proof}
    Let the BPBL pass through curves $\{\alpha_t\}_{t \in I}$, and let $B_t \subset B$ be the closed disk bounded by $\alpha_t$ for all $t \in I$. (If $\alpha_1 = \const{p}$, then define $B_1 = \{p\}$.) \Cref{lem:MinBigonRadiallyConvex} guarantees that $B$ is radially convex from $p$, and the hypotheses in this lemma imply that $B$ is convex except possibly at $p$. Thus \cref{lem:RadiallyConvexBCSP} implies that $B_t$ is also radially convex from $p$ for all $t < 1$, and so is $B_1 = \bigcap_{t < 1}B_t$.
    
    Since the BPBL only passes through simple curves, it can only converge to the constant path at $p$ or a simple geodesic loop $\gamma$ based at $p$.
    In the former case, $\{\alpha_t\}$ gives the required contraction of $\alpha$ to $p$, through curves of length at most $\length(\alpha)$. In the latter case, $\length(\gamma) \leq \length(\alpha)$. $B_1$ is convex or concave depending on its angle at $p$. If $B_1$ is concave then its interior is the bottleneck disk we seek. Otherwise it is convex, and we may apply BPFL to $\gamma$. If this process converges to a periodic geodesic, then the bottleneck disk we seek is the open disk bounded by this periodic geodesic that lies in $B_1$. Otherwise, the process converges to a point, so \cref{lem:RadiallyConvexToPoint} gives a contracting homotopy $\beta_s : \gamma \homto{\length(\alpha) + 2d} \tilde{p}$. Carrying out the homotopy $\{\alpha_t\}$ followed by $\{\beta_s\}$ gives the required contraction.
\end{proof}

\begin{figure}[h]
    \centering
    \includegraphics{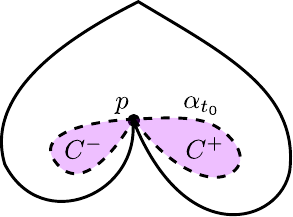}
    \caption{An illustration of the moment of self-intersection in the proof of \cref{lem:MinBigonBCSPSelfCross}.}
    \label{fig:MinBigonBCSPSelfCross}
\end{figure}

\begin{lemma}
    \label{lem:MinBigonBCSPSelfCross}
    Let $B$ be a minimizing geodesic bigon with one vertex at $p$, and with an angle of at most $\pi$ at the other vertex when measured from the inside. Parametrize $\partial B$ as a simple loop $\alpha : [0,1] \to \partial B$ based at $p$. Suppose that applying BPBL to $\alpha$ while fixing $p$ causes the resulting curves to eventually develop a self-intersection. Then either $B$ contains a bottleneck disk or $\alpha$ can be contracted to $p$ through curves in $B$ that are based at $p$ and have length at most $\length(\alpha) + 2d$.
\end{lemma}
\begin{proof}
    We will follow the general strategy used to establish Lemma~3.5 from \cite{LinearBoundsSegments}. Let $\{\alpha_t\}$ be the curves produced by BPBL applied to $\alpha$. Suppose that $\alpha_t$ develops its first self-intersection when $t = t_0$. By \cref{lem:BPBLProperties}, that self-intersection must be at $p$. Moreover, for all $t < t_0$, one of the closed disks $B_t$ bounded by $\alpha_t$ must be contained in $B$. The same lemma requires that $\bigcap_{t < t_0}B_t$ is the union of two convex closed disks $C^-$ and $C^+$ that intersect only at $p$. Let the path $\beta^\pm$ parametrize $\partial C^\pm$, and switch the signs if necessary so that $\alpha_{t_0}$ is equal to $\beta^- \concat \beta^+$ up to reparametrization. (See \cref{fig:MinBigonBCSPSelfCross}.)
    
    $B$ must be radially convex from $p$, as guaranteed by \cref{lem:MinBigonRadiallyConvex}. The hypotheses in this lemma also imply that $B$ is convex except possibly at $p$. \Cref{lem:RadiallyConvexBCSP} then requires $\bigcap_{t < t_0}B_t = C^- \cup C^+$ to also be radially convex from $p$. It can be verified that $C^-$ and $C^+$ also have to be radially convex from $p$. Applying BPFL to $\beta^\pm$ will converge either to a periodic geodesic in $C^\pm$ or a point. In the former case, the bottleneck disk we seek is the open disk bounded by this periodic geodesic that lies in $C^\pm$. In the latter case, \cref{lem:RadiallyConvexToPoint} will give a contraction $\{\beta_s^\pm\} : \beta^\pm \homto{\length(\beta^\pm) + 2d} \tilde{p}$. Assuming that $B$ does not contain a bottleneck disk, we would have both contractions $\{\beta_s^-\}$ and $\{\beta_s^+\}$. The contraction of $\alpha$ in the statement of the lemma can be carried out by applying the homotopy $\{\alpha_t\}_{t \in I}$, followed by $\{\beta^-_s \concat \beta^+\}_{s \in I}$, and then finally $\{\beta^+_s\}_{s \in I}$. The curves involved have length at most $\length(\alpha)$ in the first homotopy, at most $\length(\beta^-) + 2d + \length(\beta^+) \leq \length(\alpha) + 2d$ in the second homotopy, and at most $\length(\beta^+) + 2d \leq \length(\alpha) + 2d$ in the third homotopy. Hence the length bound in the lemma statement holds.
\end{proof}

The previous technical lemmas can now be combined to prove \cref{prop:MinBigon1ConcaveDiskOrHomotopy}.

\begin{proof}[Proof of \cref{prop:MinBigon1ConcaveDiskOrHomotopy}]
    If the angles at both vertices of $B$ are at least $\pi$ when measured from inside, then the interior of $B$ is the required bottleneck disk. Otherwise, $B$ has angle smaller than $\pi$ at one of its vertices $y$. Let $x$ be its other vertex. Reparametrize the geodesics $\gamma_i$ so that they travel from $x$ to $y$. By \cref{lem:MinBigonBCSPSimple,lem:MinBigonBCSPSelfCross}, either $B$ contains a bottleneck disk whose boundary is a geodesic loop, or there is a contraction $\partial B = \gamma_1 \concat \rev{\gamma_2} \homto{\length(\partial B) + 2d} \const{x}$ through curves in $B$. Note that $\length(\partial B) \leq \length(\gamma_1) + \length(\gamma_2) \leq 2d$. In the latter case this contraction extends to a homotopy $\{\alpha_s\} : \gamma_1 \concat \rev{\gamma_2} \concat \gamma_2 \homto{\length(\partial B) + 3d} \gamma_2$. The required homotopy is the combination of an initial homotopy $\gamma_1 \homto{3d} \gamma_1 \concat \rev{\gamma_2} \concat \gamma_2$ through curves that extend $\gamma_1$ along an arc of $\gamma_2$ and then backtracks to $y$, followed by $\{\alpha_s\}$.
\end{proof}

\section{Subadditivity of the Critical Levels}
\label{sec:Combinatorics}

Given a Riemannian sphere $M$ and points $p, q \in M$, let $v^k$ generate $H_k(\loops{pq}M;\Q)$, as explained in \cref{sec:PriorWork}. The relation in \cref{eq:PontryaginRelation_Loops} can be generalized as follows. The concatenating operation $\loops{p}M \times \loops{pq}M \xrightarrow{\concat} \loops{pq}M$ induces an algebra homomorphism $H_*(\loops{p}M) \otimes H_*(\loops{pq}M) \to H_*(\loops{pq}M)$ by the K\"unneth theorem. The isomorphism $\varphi_* : H_*(\loops{p}M) \to H_*(\loops{pq}M)$ implies that the algebra homomorphism sends $u \otimes v^k$ to $v^{k + 1}$. An argument similar to the derivation of \cref{eq:PontryaginRelation_Loops} yields the relation
\begin{equation}
    \crit(v^{k+1}) \leq \crit(v^k) + \crit(u)  \text{ for all } k \geq 0. 
    \label{eq:PontryaginRelation_kSegments1Loop}
\end{equation}

This relation is reminiscent of the notion of \emph{subadditive sequences} in Combinatorics.\footnote{For instance, see Lemma~5.1 from \cite{Ceccherini_Fekete}.} The following purely combinatorial lemma demonstrates the implications of \cref{eq:PontryaginRelation_kSegments1Loop} on the distribution of the critical levels $\crit(v^k)$.

\begin{lemma}\label{lem:SubadditiveCount}
    Let $a_0, a_1, a_2, \dotsc$ be a sequence of non-negative real numbers for which there exists a constant $A$ such that $a_{n+1} \leq a_n + A$ for all $n \geq 0$. Let the set of values attained by the sequence be $V = \{a_0, a_1, a_2, \dotsc\}$. Then for any $a_0 < r < \sup V$,
    \begin{equation}
        \Big\lvert V \cap [0,r] \Big\rvert \geq \left\lfloor\frac{r - a_0}{A}\right\rfloor + 1.
    \end{equation}
\end{lemma}

Before proving this lemma, observe that its hypotheses immediately imply that $V \cap [0,r]$ contains $a_0, a_1, \dotsc, a_{\lfloor (r - a_0)/A \rfloor}$. However, some of these values may coincide, so their contribution to the cardinality of $V \cap [0,r]$ may not suffice. In such a scenario, we would also have to consider the higher terms in the sequence that fall in the interval $[0,r]$. This corresponds to the situation where several critical levels may coincide---it may be happen that, for instance, $\lambda(v^3) = \lambda(v^5) = \lambda(v^6)$.

\begin{proof}
    Assume for the sake of contradiction that this lemma fails for some $a_0 < r < \sup V$. Then let $k = \card{V \cap [0,r]}$, where $k \leq \lfloor (r - a_0)/A\rfloor$. These $k$ distinct values divide the interval $[0,r]$ into $k+1$ subintervals. (One of these subintervals may be empty if the sequence attains the values of 0 or $r$.) $[0,a_0]$ is the union of at least one of those subintervals, and $[a_0,r]$ is the union of at most $k$ of those subintervals. We will apply the ``pigeonhole principle'' to these intervals to derive a contradiction. Let $m$ be the smallest index such that $a_m > r$, so $a_{m-1} \in V \cap [0,r]$.
    \begin{itemize}
        \item If every one of the subintervals that are contained in $[a_0,r]$ have length at most $A$, then $r - a_0 \leq kA \leq \lfloor (r - a_0)/A\rfloor A \leq r - a_0$, so equality must hold everywhere. This means that $[a_0,r]$ is the union of exactly $k$ of those subintervals, each of those $k$ subintervals having length exactly $A$. Hence those subintervals have to be $[a_0,a_0 + A], [a_0+A, a_0 + 2A], \dotsc, [r - A, r]$. Hence $V \cap [a_0,r] = \{a_0, a_0 + A, a_0 + 2A, \dotsc, r - A\}$. However, this implies that $\card{V \cap [a_0,r]} = k = \card{V \cap [0,r]}$, thus $V \cap [0,r] = \{a_0, a_0 + A, \dotsc, r - A\}$. Since $a_{m-1}$ must take one of these values, we also have $a_m \leq a_{m-1} + A \leq r$, contradicting our choice of $m$.
        
        \item Otherwise, one of those subintervals that are contained in $[a_0,r]$ has length greater than $A$. Denote this interval by $[b,c]$. Since $a_0 \leq b$, we also have $a_1 \leq b + A < c$. On the other hand, the sequence $(a_n)$ cannot take any value in $[b,c]$, so we must have $a_1 \leq b$. Similarly we can show that $a_2 \leq b$ and so on until $a_m \leq b$. However, that would contradict our assumption that $a_m > r$.
    \end{itemize}
\end{proof}

The combinatorial result above gives rise to a proof of \cref{prop:LinearBoundFirstCriticalLevel}.

\begin{proof}[Proof of \cref{prop:LinearBoundFirstCriticalLevel}]
    Fix some positive integer $k$. We will prove the bound for geodesics from $p$ to $q$, and then show how it implies the bound for geodesic loops based at $p$. The definition of a critical level implies that $\crit(v^0) = \dist{p,q}$. Let $L(k) = (k-1)\crit(u) + \dist{p,q}$. \Cref{eq:PontryaginRelation_kSegments1Loop} implies that the sequence $\crit(v^0), \crit(v^1), \crit(v^2), \dotsc$ obeys the condition in \cref{lem:SubadditiveCount} for the constant $A = \crit(u)$. Denote the set of values attained by this sequence as $V = \{\crit(v^0), \crit(v^1), \dotsc\}$. Consequently, if $L(k) < \sup V$ then applying \cref{lem:SubadditiveCount} to $r = L(k)$ gives the desired bound for geodesics from $p$ to $q$.
    
    Otherwise, $L(k) \geq \sup V$ so $V \in [0, L(k)]$. In this case, Lyusternik-Schnirelmann theory guarantees that either $\crit(v^2) < \crit(v^4) < \crit(v^6) < \dotsb$ and so on, hence $V$ is an infinite set, or if $\crit_{2i} = \crit_{2i+2}$ then $\loops{pq}M$ contains an infinite sequence of geodesics from $p$ to $q$ whose lengths approach that value. Either way, we obtain the bound for geodesics from $p$ to $q$.
    
    The bound for geodesic loops based at $p$ follows from the same arguments applied to the scenario where $p = q$, as long as we compensate for the fact that the shortest one of these geodesic loops would actually be $\const{p}$, by raising the length bound to $L(k+1)$.
\end{proof}

\section{Morse Length Functionals and the Proof of Theorem~\ref{thm:TighterLinearBound}}
\label{sec:MorseHomology}

The final scenario mentioned in \cref{thm:TighterLinearBound} involves a Morse length functional, which is addressed in the following lemma.

\begin{lemma}
    \label{lem:LinearBound_MorseCobordism}
    Let $M$ be a Riemannian sphere, and fix two points $p, q \in M$. Assume that the metric on $M$ is such that the geodesics from $p$ to $q$ are nondegenerate. Then $M$ has at least $k$ distinct geodesics from $p$ to $q$ of length at most $\lfloor \frac{k-1}2\rfloor\crit(u^2) + (k - 1 \bmod2)\crit(u) + \dist{p,q}$, where $u$ is a generator of $H_1(\loops{p}M;\Q)$.
    
    In addition, $M$ also has at least $k$ distinct non-constant geodesic loops based at $p$ of length at most $\lfloor k/2\rfloor\crit(u^2) + (k \bmod2)\crit(u)$.
\end{lemma}

\begin{proof}
    We will prove this theorem for geodesics from $p$ to $q$ in the case where $k$ is even, because the odd case requires only a minor modification to the argument. The bound for geodesic loops based at $p$ will then follow by setting $p = q$, while increasing the value of $k$ in the bound by 1 to compensate for the fact that one of the geodesics from $p$ to itself is the constant loop. The length bound in this situation is $L = (k/2 - 1)\crit(u^2) + \crit(u) + \dist{p,q}$.     Define the homology classes $v^i \in H_i(\loops{pq}M;\Q)$ as in \cref{sec:PriorWork}. It follows from \cref{eq:PontryaginRelation_kSegments1Loop,eq:CriticalLevelsEvenBoundU2} that for any integer $0 
    \leq i \leq k/2-1$,
    \begin{align}
        \crit(v^{2i}) &\leq i\crit(u^2) + \dist{p,q} \leq L\\
        \crit(v^{2i+1}) &\leq \crit(v^{2i}) + \crit(u)
        \leq i\crit(u^2) + \crit(u) + \dist{p,q} \leq L,
    \end{align}
    therefore $\crit(v^0), \crit(v^1), \dotsc, \crit(v^{k-1}) \leq L$. Since the length functional on $\loops{p}M$ is Morse by hypothesis, a standard argument implies that the geodesic $\gamma^i$ obtained from $v^i$ by the Birkhoff minimax principle has Morse index $i$. Hence we have $k$ distinct geodesics $\gamma_0, \gamma_1, \dotsc, \gamma_{k-1}$ of length at most $L$.
\end{proof}

Finally, we can tie all of our previous results together to prove \cref{thm:TighterLinearBound}.

\begin{proof}[Proof of \cref{thm:TighterLinearBound}]
    If $p \neq q$, then we apply the $n = 4$ case of \cref{thm:SweepoutOrBottleneckDisks}: either $M \setminus \{p\}$ contains 4 disjoint bottleneck disks whose boundaries may consist of at most four geodesics, or there is a sweepout of $M$ whose longest induced meridian has length at most $7d + o(\epsilon)$ and whose shortest induced meridian has length at most $d + o(\epsilon)$. In the former case, at least 2 of the disks would not contain $q$, so their complement would be a convex annulus---or some convex set that is homotopy equivalent to an annulus---containing $p$ and $q$. Thus the argument in \cref{sec:PriorWork} (whose details can be found in section 2 of \cite{LinearBoundsSegments}) gives enough distinct geodesics from $p$ to $q$ ``winding around'' the convex annulus to satisfy the bound in the theorem. In the latter case, \cref{lem:MeridiansBoundCriticalLevels} guarantees that $\crit(u) \leq 8d + o(\epsilon)$ and $\crit(u^2) \leq 14d + o(\epsilon)$, therefore, $\crit(u) \leq 8d$ and $\crit(u^2) \leq 14d$. \Cref{prop:LinearBoundFirstCriticalLevel} then implies that there are at least $k$ geodesics from $p$ to $q$ of length at most $8(k-1)d + \dist{p,q} \leq (8k-7)d \leq 8kd$.
    
    Choose an integer $2 \leq m \leq \infty$. Now we tighten the above bound for a set of metrics that is generic in the $C^m$ topology. It was proven in \cite{Bettiol_GenericNondegen} that for such a generic set of metrics $g$, the manifold $M = (\Sp,g)$ has no degenerate geodesics from $p$ to $q$.  \Cref{lem:LinearBound_MorseCobordism} then implies that $M$ has $k$ distinct geodesic from $p$ to $q$ of length at most $ (7k - 6 - (k\bmod2))d + \dist{p,q}$.
    
    The $p = q$ case of this theorem can be proven in nearly the same manner as above, except that we apply the $n = 2$ case of \cref{thm:SweepoutOrBottleneckDisks} instead.
\end{proof}

\bibliographystyle{amsplain}

\begin{thebibliography}{10}

\bibitem{Ballmann_LyusternikSchnirelmann}
W.~Ballmann, G.~Thorbergsson, and W.~Ziller, \emph{Existence of closed
  geodesics on positively curved manifolds}, Journal of differential geometry
  \textbf{18} (1983), no.~2, 221--252.

\bibitem{Beach_FreeLoopSpace}
I.~Beach, A.~Nabutovsky, and R.~Rotman, \emph{Quantitative {M}orse theory on
  free loop spaces on {R}iemannian 2-spheres}, In preparation.

\bibitem{Bettiol_GenericNondegen}
Renato~G Bettiol and R.~Giamb{\`o}, \emph{Genericity of nondegenerate geodesics
  with general boundary conditions}, Topological Methods in Nonlinear Analysis
  \textbf{35} (2010), no.~2, 339--365.

\bibitem{BirkhoffHestenes_Minimax}
G.~D. Birkhoff and M.~R. Hestenes, \emph{Generalized minimax principle in the
  calculus of variations}, Proceedings of the National Academy of Sciences of
  the United States of America \textbf{21} (1935), no.~2, 96.

\bibitem{Buchner_CutLocus}
M.~A. Buchner, \emph{The structure of the cut locus in dimension less than or
  equal to six}, Compositio Mathematica \textbf{37} (1978), no.~1, 103--119.

\bibitem{Ceccherini_Fekete}
T.~Ceccherini-Silberstein, M.~Coornaert, and F.~Krieger, \emph{An analogue of
  {F}ekete's lemma for subadditive functions on cancellative amenable
  semigroups}, Journal d'Analyse Math{\'e}matique \textbf{124} (2014), no.~1,
  59--81.

\bibitem{Cheeger_InjectivityRadius}
J.~Cheeger, \emph{Finiteness theorems for {R}iemannian manifolds}, American
  Journal of Mathematics \textbf{92} (1970), no.~1, 61--74.

\bibitem{ColdingDeLellis_Minmax}
T.~H. Colding and C.~De~Lellis, \emph{The min-max construction of minimal
  surfaces}, Surveys in Differential Geometry \textbf{8} (2003), 75--107.

\bibitem{Croke}
C.~B. Croke, \emph{Area and the length of the shortest closed geodesic},
  Journal of Differential Geometry \textbf{27} (1988), no.~1, 1--21.

\bibitem{Felix_RationalHomotopyTheory}
Y.~F{\'e}lix, S.~Halperin, and J.-C. Thomas, \emph{Rational homotopy theory},
  vol. 205, Springer Science \& Business Media, 2012.

\bibitem{FrankelKatz}
S.~Frankel and M.~Katz, \emph{The {M}orse landscape of a {R}iemannian disk},
  Annales de l'institut Fourier, vol.~43, 1993, pp.~503--507.

\bibitem{Liokumovich_LongSweepouts}
Y.~Liokumovich, \emph{Spheres of small diameter with long sweep-outs},
  Proceedings of the American Mathematical Society \textbf{141} (2013), no.~1,
  309--312.

\bibitem{LNR_PeriodicGeodesicsS2}
Y.~Liokumovich, A.~Nabutovsky, and R.~Rotman, \emph{Lengths of three simple
  periodic geodesics on a {R}iemannian 2-sphere}, Mathematische Annalen
  \textbf{367} (2017), no.~1-2, 831--855.

\bibitem{Lyusternik_Book}
L.~Lyusternik, \emph{The topology of function spaces and the calculus of
  variations in the large}, vol.~16, American Mathematical Soc., 1967.

\bibitem{LyusternikSchnirelmann_ThreeGeodesics}
L.~Lyusternik and L.~Schnirelmann, \emph{Sur le probl{\`e}me de trois
  g{\'e}od{\'e}siques ferm{\'e}es sur les surfaces de genre 0}, CR Acad. Sci.
  Paris \textbf{189} (1929), 269--271.

\bibitem{MarquesNeves_MinimalVarieties}
F.~C. Marques and A.~Neves, \emph{Topology of the space of cycles and existence
  of minimal varieties}, Surveys in differential geometry \textbf{21} (2016),
  no.~1, 165--177.

\bibitem{MorseTheory_Milnor}
J.~Milnor, \emph{{M}orse theory. ({AM}-51)}, vol.~51, Princeton University
  Press, 1963.

\bibitem{Myers_AnalyticCutLocus}
S.~B. Myers, \emph{Connections between differential geometry and topology},
  Proceedings of the National Academy of Sciences of the United States of
  America \textbf{21} (1935), no.~4, 225.

\bibitem{LinearBounds}
A.~Nabutovsky and R.~Rotman, \emph{Linear bounds for lengths of geodesic loops
  on {R}iemannian 2-spheres}, Journal of Differential Geometry \textbf{89}
  (2011), no.~2, 217--232.

\bibitem{Nabutovsky_QuantitativeMorse}
\bysame, \emph{Length of geodesics and quantitative {M}orse theory on loop
  spaces}, Geometric and Functional Analysis \textbf{23} (2013), no.~1,
  367--414.

\bibitem{LinearBoundsSegments}
\bysame, \emph{Linear bounds for lengths of geodesic segments on {R}iemannian
  2-spheres}, Journal of Topology and Analysis \textbf{5} (2013), no.~04,
  409--438.

\bibitem{Schwartz}
A.~S. Schwartz, \emph{Geodesic arcs on {R}iemannian manifolds}, Uspekhi
  Matematicheskikh Nauk \textbf{13} (1958), no.~6, 181--184.

\bibitem{Serre}
J.-P. Serre, \emph{Homologie singuli{\`e}re des espaces fibr{\'e}s}, Annals of
  Mathematics (1951), 425--505.

\end{thebibliography}
\providecommand{\bysame}{\leavevmode\hbox to3em{\hrulefill}\thinspace}
\providecommand{\MR}{\relax\ifhmode\unskip\space\fi MR }
\providecommand{\MRhref}[2]{%
  \href{http://www.ams.org/mathscinet-getitem?mr=#1}{#2}
}
\providecommand{\href}[2]{#2}

\end{document}